\documentclass[11pt,a4paper,reqno]{amsart}

\usepackage{amsmath}
\usepackage{amsthm}
\usepackage{amsfonts, amssymb}
\usepackage[all]{xy}
\usepackage{url}
\usepackage{algorithm}
\usepackage[end]{algpseudocode}
\usepackage{framed}
\usepackage{faktor}
\usepackage{tabularx}

\usepackage{latexsym}

\usepackage{graphicx}
\usepackage{graphics}
\usepackage[scriptsize,bf]{caption}
\usepackage{float}
\usepackage{enumerate}
\usepackage{verbatim}
\usepackage{color}
\usepackage{subfig}
\usepackage{tikz}

\theoremstyle{plain}
\newtheorem{lema}{Lemma}[section]
\newtheorem{prop}[lema]{Proposition}
\newtheorem{teo}[lema]{Theorem}
\newtheorem{conj}[lema]{Conjecture}

\newtheorem{coro}[lema]{Corollary}

\theoremstyle{definition}

\newtheorem{preg}[lema]{Question}

\newtheorem{obs}[lema]{Remark}

\makeatletter
\@namedef{subjclassname@2020}{%
  \textup{2020} Mathematics Subject Classification}
  
\newcommand{\multiline}[1]{%
  \begin{tabularx}{\dimexpr\linewidth-\ALG@thistlm}[t]{@{}X@{}}
    #1
  \end{tabularx}
}  
  
\makeatother

\begin{document}

\title[On sequences arising from randomizing subtraction games]{On sequences arising from randomizing \\subtraction games}

\author[N. A. Capitelli]{Nicol\'as Capitelli}
\author[F. Somma]{Francisco Somma\\
\textit{\scriptsize
U\MakeLowercase{niversidad} N\MakeLowercase{acional de} L\MakeLowercase{uj\'an}, D\MakeLowercase{epartamento de} C\MakeLowercase{iencias} B\MakeLowercase{\'asicas}, A\MakeLowercase{rgentina.}}}

\thanks{\textit{Corresponding address:} {\color{blue}ncapitelli@unlu.edu.ar}}

\subjclass[2020]{91A46, 91A60, 11B83, 40A05}

\keywords{Sequences, Combinatorial games, Subtraction games}

\thanks{\textit{This research was partially supported by the Department of Basic Sciences, UNLu (Disposition CDD-CB 148/18) and University of Luj\'an (Resolutions REC 224/19 and PHCS 084/20)}}

\begin{abstract} In this article, we study the behaviour of a broad family of real sequences derived from randomized one-pile subtraction games. For any subtraction set $S$, we allow any valid number of chips $s\in S$ to be removed at equal probability at any given position and we study the sequences $(a_n^S)_{n\in\mathbb{N}}$ representing the probability of winning the game from a position with $n$ chips. We characterize these sequences in terms of linear recurrence relations and examine their behavior as $n\rightarrow\infty$ for all finite $S$. We fully solve the cases for subtraction sets of fewer than 3 elements and partially complete the general case for arbitrary $S$.\end{abstract}

\maketitle

A (one-pile) subtraction game with subtraction set $S\subset\mathbb{N}$ is a combinatorial game which is played as follows. From a pile with $n$ chips, two players alternate moves. A move consists of removing $s$ chips from the pile where $s \in S$. The player that has no moves  left available loses. Subtraction games are fundamental to combinatorial game theory and have a rich theoretical background. However, analyzing them for generic subtraction sets can be quite challenging, with many open problems still remaining (see e.g. \cite{Nowa}). For a comprehensive understanding of combinatorial games, we refer the reader to the seminal work \cite{WW} (see also \cite[\S 1]{Ferguson})

In this article, we examine the behaviour of a large family of real sequences which we derive by randomizing subtraction games. Specifically, we allow one such game to be played stochastically, with an equal probability of removing any valid number of chips from the subtraction set $S$ at a given position. We then let $a_n^S$ be the probability that the next player to move will win this game from a position with a pile of $n$ chips. For each finite $S$, we characterize the sequences $(a_n^S)_{n\in\mathbb{N}}$ in terms of linear recurrence relations and examine their behavior as $n\rightarrow\infty$. For the case $|S|= 2$, we univoquely characterize their convergence and we explicit exhibit the limit value of the subsequences of even and odd terms when $\lim_{n\rightarrow\infty}(a_n^S)_{n\in\mathbb{N}}$ does not exist (Theorem \ref{Thm:Main2}). As a general result, we prove convergence for all finite subtraction sets $S$ containing at least an even number (Theorem \ref{Thm:MainGeneral}).

\section{A family of sequences}

Let $S\subset\mathbb{N}$ be a non-empty set. We randomize the one-pile subtraction game of $n$ chips with subtraction set $S$ by playing it as follows.
\begin{enumerate}
\item First player to move starts from the starting position (an $n$-chips pile).
\item Players alternate moves.
\item At a given position (number of chips in the pile) $0\leq m\leq n$, an element $k'\in \{k\in S\,|\, k\leq m\}$ is randomly selected with equal probability, and the next player to move removes $k'$ chips from the pile.
\item The first player who is confronted with a position $m$ for which $\{k\in S\,|\, k\leq m\}=\emptyset$, loses.
\end{enumerate}

These randomly-played games are a particular case of those studied in \cite{CapitelliPrivitelli}, which involve playing through noisy channels that can introduce errors in transmission.

Our main concern is the sequence $a_n^S$ which represents the probability of winning for the next player $N$ to move in a game of $n$ chips and subtraction set $S$. Note that, for $n\geq \max\{S\}$, $a_n^S$ can be computed recursively as follows.
\begin{eqnarray}\label{Eq:Recursive}a_n^S&=&\displaystyle\sum_{i\in S}\ Prob(\text{$N$ wins from $n$ $|$ $i$ is selected}) \cdot Prob(\text{$i$ is selected}) \\
&=&\displaystyle\sum_{i\in S}\ (1-Prob(\text{$P$ wins from $n-i$})) \cdot \frac{1}{|S|}\nonumber  \\
&=&\displaystyle\sum_{i\in S} (1-a_{n-i}^S)\cdot \frac{1}{|S|} \nonumber\\
&=&1-\frac{1}{|S|}\displaystyle\sum_{i\in S} a_{n-i}^S \nonumber\end{eqnarray}
Here $P$ denotes the previous player, $Prob(X)$  the probability of event $X$ and $Prob(X|Y)$ the conditional probability of $X$ given $Y$. Once more, this line of reasoning mirrors that found in \cite{CapitelliPrivitelli}.

We see that the sequence $a_n^S$ satisfies a non-homogeneous linear recurrence of order $\max\{S\}$. To determine the initial values of $a_i^S$ for $i < \max\{S\}$, we proceed as follows. Let $S = \{k_1, \ldots, k_t\}$ with $k_1 < \cdots < k_t$. Thus, for $i < k_t$, $a_i^S$ represents the probability of the next player winning from position $i$ in a game with subtraction set $S - \{k_t\} = \{k_1, \ldots, k_{t-1}\}$. Therefore, the sequence $a_n^S$ can be determined inductively based on the size of $S$. For $t = 1$, there is no randomness, and $a_n^S$ coincides with the Sprague-Grundy function of the subtraction game with subtraction set $S = {k_1}$ (see \eqref{Eq:Case_t=1} below for an explicit formula for this function). For $t>1$ we let:
\begin{equation}\label{Eq:Main}
    a_n^S=\begin{cases}a_n^{S-\{k_t\}}& 0\leq n< k_t\\1-\frac{1}{t}\sum_{i=1}^t a_{n-k_i}^S & n\geq k_t \end{cases}
\end{equation}

Our main focus in this article is to analyze the behavior of $a_n^S$ for large $n$, which serves as an indicator of the game's fairness. It's straighfoward to verify that convergence of $a_n^S$ implies convergence to $\frac{1}{2}$. However, this will also emerge as a consequence of our analysis of the recursive relation \eqref{Eq:Main}. A solution for this non-homogenous linear recurrence is of the form \begin{equation}\label{Eq:ClosedForm}a_n^S=\alpha_1 z_1^n+\alpha_2 z_2^n+\cdots + \alpha_l z_{k_t}^n+\frac{1}{2}\end{equation} where $\{z_i\}_{1\leq i\leq k_t}$ are the complex roots of the characteristic polynomial \begin{equation}\label{Eq:Characteristic}\chi_S(x)=tx^{k_t}+\sum_{i=1}^{t} x^{k_t-k_i}\end{equation} 
and $\alpha_i\in\mathbb{C}$ are the solutions to the system
\begin{equation}\label{Eq:System}\begin{cases}
    \alpha_1+\alpha_2+\cdots+\alpha_{k_t}+\frac{1}{2}=a_0^S\\
\alpha_1z_1+\alpha_2z_2+\cdots+\alpha_{k_t}z_{k_t}+\frac{1}{2}=a_1^S\\
\alpha_1z_2^1+\alpha_2z_2^2+\cdots+\alpha_{k_t}z_{k_t}^2+\frac{1}{2}=a_2^S\\
\hspace{1.2in}\vdots\\
\alpha_1z_1^{{k_t}-1}+\alpha_2z_2^{{k_t}-1}+\cdots+\alpha_{k_t}z_{k_t}^{{k_t}-1}+\frac{1}{2}=a_{{k_t}-1}^S
\end{cases}\end{equation}
(see, e.g., \cite[\S 1]{Knuth} for a glimpse of recursive sequences theory). In the following sections we'll study the behaviour of $(a_n^S)_{n\in\mathbb{N}_0}$, completely solving the case $t=2$ and partially the general case for finite $S$.

\begin{obs}\label{Lema:Converges_1/2} Note that in order to study the limit of $(a_n^S)_{n\in\mathbb{N}_0}$, we may restrict our attention to the case where $\gcd\{k_1,\ldots,k_t\}=1$. Indeed, let $m=\gcd\{k_1,\ldots,k_t\}$. Consider the subsequence defined as $b_n=a^S_{mn}$. Then $$b_n=a^S_{mn}=1-\frac{1}{t}\sum_{i=1}^t a_{mn-k_i}^S$$
    Since $mn-k_i=mn-m(k_i/m)=m(n-k_i/m)$ then $a_{mn-k_i}^S=b_{n-k_i/m}$. Hence, $$b_n=1-\frac{1}{t}\sum_{i=1}^t b_{n-k_i/m},$$ which is the sequence associated to the one-pile subtraction game of $n$ chips and subtraction set $S'=\{\frac{k_1}{m},\ldots,\frac{k_t}{m}\}$. Since $\lim_{n\rightarrow \infty} a_n=\lim_{n\rightarrow \infty} b_n$ the claim is proven. 
\end{obs}

As it was pointed out earlier, for $S=\{k\}$, the game is devoid of randomness, with either no moves or just one available. In this case, the sequence adopts the following straightforward closed form:
\begin{equation}\label{Eq:Case_t=1}a_n^{\{k\}}=\dfrac{1+(-1)^{\lfloor \frac{n}{k}\rfloor+1}}{2}\end{equation}
That is, $a_n^{\{k\}}$ consists of alternating blocks of zeros and ones, each block being of size $k$.

\section{The case $S=\{k,l\}$}
In this section with shall completely solve the case of a subtraction set with two elements. Namely, we will univoquely determine for which subtraction sets $S$ the sequence $a_n^S$ converges and for which it does not. Moreover, in the latter scenario, the subsequences of even and odd terms converge, and we explicitly exhibit these limits for each $S$. Our approach to solving this case mirrors our strategy for addressing the general case.

Throughout this section we let $S=\{k,l\}\subset\mathbb{N}$ with $k<l$ and $\gcd\{k,l\}=1$ (see Remark \ref{Lema:Converges_1/2}). We also shall write $a_n$ for $a_n^S=a_n^{\{k,l\}}$ for simplicity. In this case, Equations \eqref{Eq:Main} and \eqref{Eq:Case_t=1} give us the following recursive sequence.
\begin{equation}\label{Eq:Main_kl}
    a_n=\begin{cases}\frac{1+(-1)^{\lfloor \frac{n}{k}\rfloor+1}}{2}& 0\leq n< l\\1-\frac{1}{2}(a_{n-k}+a_{n-l}) & n\geq l \end{cases}
\end{equation}
The following properties of the characteristic polynomial $\chi(z)=2z^l+z^{l-k}+1$ of this recurrence will help us characterize the convergence of $(a_n)_{n\in\mathbb{N}_0}$ for each election of $k$ and $l$. In what follows we let $G_i\subset\mathbb{C}$ denote the group of $i$th roots of unity.

\begin{prop}\label{Prop:Roots2}
    All the complex roots of $\chi(x)$ have norm less than 1, with the exception of the root $z=-1$ which happens only for the case $k,l\equiv 1 \, (2)$. Furthermore, in this last case, all roots are simple.

\begin{proof}
    Let $z$ be a root of $\chi(x)$. By Lagrange's bound (see e.g. \cite[IV]{Yap}) we have that \begin{equation}|z|\leq \max\left\{1,\sum_{i=0}^{d-1}\left|\frac{c_i}{c_l}\right|\right\}=1,\end{equation} where $c_0,\ldots,c_{d_t}$ denote the coefficients of $\chi(x)$. Suppose $|z|=1$. Write $-1=2z^{l}+z^{l-k}=z^{l-k}(2z^{k}+1)$. Since $|z|=1$ then $|2z^{k}+1|=1$. Therefore $2z^{k}+1$ lies in the intersection of the circle of center $(1,0)$ and radius $2$ and the circle of center $(0,0)$ and radius 1. Hence $2z^{k}+1=-1$. We conclude that $z^{l-k}=1$ and that $z^{l}=-1$. Particularly, $z$ is a root of unity of order $2l$. Therefore $z\in G_{l-k}\cap G_{2l}=G_{\gcd\{l-k,2l\}}$. Since $\gcd\{l-k,l\}=1$ then $$G_{\gcd\{l-k,2l\}}=\begin{cases}
        G_1 & l-k\equiv 1\, (2)\\
        G_2 & l-k\equiv 0\, (2)
    \end{cases}$$
    Now, $\chi(1)\neq 0$ so if $l-k\equiv 1\, (2)$ the $\chi$ has no roots of norm 1 and if $l-k\equiv 0\, (2)$ the only possible root of norm 1 is $z=-1$. Since $\chi(-1)=0$ for this last case (which implies $l$ odd), the result is settled.

    It remains to show the simplicity of the all of the roots for the case $l$ and $k$ odd. In this case, $z=-1$ is a simple root since  $\chi'(-1)\neq 0$. Let $z\neq -1$ and suppose $z$ is a root of $\chi'(x)$. Then $z^{l-k-1}(2lz^{k}+l-k)=0$, which implies that $z^{k}\in\mathbb{R}$. If $z$ is also a root of $\chi(x)$ then $0=2z^l+z^{l-k}+1=z^{l-k}(2z^{k}+1)+1$ which implies $z^{l-k}\in\mathbb{R}$ and $z^l\in\mathbb{R}$.
This means that $\arg(z)=\frac{r\pi}{l-k}$ for some $0\leq r< 2(l-k)$ and $\arg(z)=\frac{r'\pi}{l}$ for some $0\leq r'< 2l$ simultaneously. Hence $lr=(l-k)r'$. Since $\gcd\{l-k,l\}=1$ then $l-k$ divides $r$. Since $0\leq r< 2(l-k)$ then either $r=0$ o $r=l-k$, from which $\arg(z)=0$ or $\arg(z)=\pi$. So $z\in\mathbb{R}$. In the first case, $z>0$ but $\chi(x)$ doesn't have real positive roots. In the second case, since $z\neq -1$ and $-1$ is a root for this case, then it must cancel the polynomial $\chi(x)/(x+1)$. An easy verification shows that \begin{equation}\label{Eq:Division}\frac{\chi(x)}{x+1}=\sum_{i=0}^{l-k-1} (-1)^i x^i + 2\sum_{i=l-k}^{l-1} (-1)^i x^i\end{equation}

But since $z<0$ then $(-1)^i z^i = (-z)^i >0$ for all $i$. Therefore, $z$ does not anulate this quotient and cannot be a multiple root of $\chi(x)$.\end{proof}
\end{prop}

\begin{coro}\label{Coro:1stCoro}
    The sequence $(a_n)_{n\in\mathbb{N}_0}$ defined in \eqref{Eq:Main_kl} converges to $\frac{1}{2}$ if $l\not\equiv k\, (2)$.

    \begin{proof}
        If $l\not\equiv k (2)$ then Proposition \ref{Prop:Roots2} implies that all roots of $\chi(z)=2z^l+z^{l-k}+1$ have norm less than $1$ (since if $l$ is odd then $l-k$ is also odd). Hence, $\lim_{n\rightarrow \infty}a_n=\frac{1}{2}$ by equation \eqref{Eq:ClosedForm}.
    \end{proof}
\end{coro}

The case $l,k\equiv 1 \, (2)$ remains to be addressed. We tackle this in the following steps. Under these conditions, Proposition \ref{Prop:Roots2} gives $z_1=-1$ as the only  root of $\chi(x)$ of norm 1, while all remaining roots possess norm less than 1. From Equation \eqref{Eq:ClosedForm} we conclude that \begin{equation}\label{Eq:alpha_1+1/2}\lim_{n\rightarrow\infty} a_n = \lim_{n\rightarrow\infty} \alpha_1 (-1)^n+\frac{1}{2}.\end{equation} We shall explicity compute $\alpha_1$ next. Note that the linear system \eqref{Eq:System} becomes:

\begin{equation}\label{Eq:MatrixForm}\begin{pmatrix}
    1&1&1&\cdots&1\\
    -1&z_2&z_3&\cdots&z_l\\
    1&z_2^2&z_3^2&\cdots&z_l^2\\
    \vdots& \vdots& \vdots& & \vdots\\
    -1&z_2^{l-2}&z_3^{l-2}&\cdots&z_l^{l-2}\\
    1&z_2^{l-1}&z_3^{l-1}&\cdots&z_l^{l-1}
\end{pmatrix}\begin{pmatrix}
    \alpha_1\\
    \alpha_2\\
    \alpha_3\\
    \vdots\\
    \alpha_{l-1}\\
    \alpha_l
\end{pmatrix}=\begin{pmatrix}
    a_0-\frac{1}{2}\\
    a_1-\frac{1}{2}\\
    a_2-\frac{1}{2}\\
    \vdots\\
    a_{l-2}-\frac{1}{2}\\
    a_{l-1}-\frac{1}{2}
\end{pmatrix}\end{equation}

Since all roots of $\chi(x)$ are simple (again by Proposition \ref{Prop:Roots2}) then this Vandermonde matrix is invertible. The inverse is given by

$$\begin{pmatrix}
    L_{00}&\cdots&L_{l-1 0}\\
    \vdots& & \vdots\\
    L_{0 l-1}&\cdots&L_{l-1 l-1}
\end{pmatrix}$$
where $L_j(x)=\sum_{i=0}^{l-1}L_{ij}x^i=\frac{f(x)}{(x-z_{j+1})f'(z_{j+1})}$ and $f(x)=\prod_{i=1}^{l}(x-z_i)$ (see e.g. \cite[\S0.9]{HornJohnson} for information on Vandermonde matrices). We conclude that:
$$\begin{pmatrix}
    \alpha_1\\
    \vdots\\
    \alpha_l
\end{pmatrix}=\begin{pmatrix}
    L_{00}&\cdots&L_{l-1 0}\\
    \vdots& & \vdots\\
    L_{0 l-1}&\cdots&L_{l-1 l-1}
\end{pmatrix}\begin{pmatrix}
    a_0-\frac{1}{2}\\
    \vdots\\
    a_{l-1}-\frac{1}{2}
\end{pmatrix}$$
Since we are solely interested in $\alpha_1$ we only need consider \begin{equation}\label{L_0}L_0(x)=\dfrac{f(x)}{(x-z_1)f'(z_1)}=\dfrac{f(x)}{(x+1)f'(-1)}.\end{equation} 
Since $\chi(x)=2f(x)$ then, by equation \eqref{Eq:Division}, we get
$$\dfrac{f(x)}{(x+1)}=\frac{1}{2}\sum_{i=0}^{l-k-1} (-1)^i x^i + \sum_{i=l-k}^{l-1} (-1)^i x^i$$
Also, $f'(z_1)=f'(-1)=\frac{1}{2}(l+k)=\frac{l+k}{2}$. Hence $$L_0(x)=\frac{2}{l+k}\left(\frac{1}{2}\sum_{i=0}^{l-k-1} (-1)^i x^i + \sum_{i=l-k}^{l-1} (-1)^i x^i\right)$$

We conclude that the coefficients of $L_0(x)$ are then given by $$L_{i0}=\begin{cases}
    \dfrac{(-1)^i}{l+k} & 0\leq i<l-k\\
    &\\
    \dfrac{2(-1)^i}{l+k} & l-k\leq i< l\\
\end{cases}$$

Hence \begin{eqnarray}\alpha_1&=&\sum_{i=0}^{l-k-1}\dfrac{(-1)^i}{l+k}\left(a_i-\frac{1}{2}\right)+\sum_{i=l-k}^{l-1}\dfrac{2(-1)^i}{l+k}\left(a_i-\frac{1}{2}\right)\nonumber\\
&=&\dfrac{1}{l+k}\sum_{i=0}^{l-k-1}(-1)^i\left(a_i-\frac{1}{2}\right)+\dfrac{2}{l+k}\sum_{i=l-k}^{l-1}(-1)^i\left(a_i-\frac{1}{2}\right)\nonumber\\
&=&\dfrac{1}{l+k}\left(\sum_{i=0}^{l-k-1}(-1)^i\left(a_i-\frac{1}{2}\right)+2\sum_{i=l-k}^{l-1}(-1)^i\left(a_i-\frac{1}{2}\right)\right)\nonumber\\
&=&\dfrac{1}{l+k}\left(\sum_{i=0}^{l-k-1}(-1)^ia_i-\sum_{i=0}^{l-k-1}(-1)^i\frac{1}{2}+2\sum_{i=l-k}^{l-1}(-1)^ia_i-2\sum_{i=l-k}^{l-1}(-1)^i\frac{1}{2}\right)\nonumber\\
&=&\dfrac{1}{l+k}\left(\sum_{i=0}^{l-k-1}(-1)^ia_i+2\sum_{i=l-k}^{l-1}(-1)^ia_i-\sum_{i=0}^{l-k-1}(-1)^i\frac{1}{2}-2\sum_{i=l-k}^{l-1}(-1)^i\frac{1}{2}\right)\nonumber\end{eqnarray}
Now, since $l-k-1$ is odd then the sum $\sum_{i=0}^{l-k-1}(-1)^i\frac{1}{2}$ is zero and since $l-k$ and $l-1$ are even then the sum $\sum_{i=l-k}^{l-1}(-1)^i\frac{1}{2}$ is $\frac{1}{2}$. We conclude
$$\alpha_1=\dfrac{1}{l+k}\left(\sum_{i=0}^{l-k-1}(-1)^ia_i+2\sum_{i=l-k}^{l-1}(-1)^ia_i-1\right)=\dfrac{1}{l+k}\left(\sum_{i=0}^{l-1}(-1)^ia_i+\sum_{i=l-k}^{l-1}(-1)^ia_i-1\right)$$
We compute these sums in the following 

\begin{lema} Let $a_j=\dfrac{1+(-1)^{\lfloor \frac{j}{k}\rfloor+1}}{2}$ for $0\leq j\leq l-1$. Then,

\begin{enumerate}
    \item $\displaystyle\sum_{i=0}^{l-1}(-1)^ia_i=\begin{cases}
        -\frac{1}{2}\left\lfloor\frac{l}{k}\right\rfloor & \left\lfloor\frac{l}{k}\right\rfloor\equiv 0 \, (2)\\ 
        -\frac{1}{2}\left(\left\lfloor\frac{l}{k}\right\rfloor -1\right) & \left\lfloor\frac{l}{k}\right\rfloor\equiv 1 \, (2)
    \end{cases}$ \vspace{0.1in}
    \item $\displaystyle\sum_{i=l-k}^{l-1}(-1)^ia_i=\begin{cases}
        -k\left(\left\lfloor\frac{l}{k}\right\rfloor +1 \right)+l & \left\lfloor\frac{l}{k}\right\rfloor\equiv 0 \, (2)\\ 
        -l+k\left\lfloor\frac{l}{k}\right\rfloor  & \left\lfloor\frac{l}{k}\right\rfloor\equiv 1 \, (2)
    \end{cases}$
\end{enumerate}
    \begin{proof}
        Let $s$ be such that $sk\leq i\leq (s+1)k-1$ and note that $$(-1)^ia_i=\begin{cases}
            0 & s\equiv 0 \, (2)\\
            (-1)^i & s\equiv 1 \, (2)
        \end{cases}$$
        So, we have blocks of zeros and $\pm 1$ alternating. All the blocks of $\pm 1$ begin at $i=sk$ with odd $s$, so the leading term of each of these blocks is $-1$. Since there are $(s+1)k-1-sk+1=k$ integers in each interval and $k$ is odd then the sum of each block of $\pm 1$ is $-1$.
        Let us write $l=k\lfloor\frac{l}{k}\rfloor+r$, where $r$ is the remainder of the division of $l$ by $k$. Then, we have $\lfloor\frac{l}{k}\rfloor$ alternating blocks of zeros and $\pm 1$ and a truncated block of size $r$. More concretely,
\begin{itemize}
    \item If $\lfloor\frac{l}{k}\rfloor\equiv 0 \, (2)$ then the truncated block is of zeros and we have $\frac{1}{2}\lfloor\frac{l}{k}\rfloor$ complete blocks of $\pm 1$. In this case, $\sum_{i=0}^{l-1}(-1)^ia_i=-\frac{1}{2}\lfloor\frac{l}{k}\rfloor$.
    \item If $\lfloor\frac{l}{k}\rfloor\equiv 1 \, (2)$ then the truncated block is of $\pm 1$ and we have $\frac{1}{2}(\lfloor\frac{l}{k}\rfloor-1)$ complete blocks of $\pm 1$. The truncated block has $r=l-k\lfloor\frac{l}{k}\rfloor$ $\pm 1$. Since $l$, $k$ and $\lfloor\frac{l}{k}\rfloor$ are odd then $r$ is even. So the sum of the truncated block of $r$ $\pm 1$ cancels and $\sum_{i=0}^{l-1}(-1)^ia_i=-\frac{1}{2}(\lfloor\frac{l}{k}\rfloor-1)$.
\end{itemize}
Now, for the sum $\sum_{i=l-k}^{l-1}(-1)^ia_i$ note that there are $l-1-(l-k)+1=k$ sumands. Since if $\lfloor\frac{l}{k}\rfloor\equiv 0 \, (2)$ the truncated block is of zeros then we must have $r$ zeros and $k-r$ $\pm 1$. Since $r=l-k\lfloor\frac{l}{k}\rfloor$ then we have $k-r=k-(l-k\lfloor\frac{l}{k}\rfloor)=k(\lfloor\frac{l}{k}\rfloor+1)-l$ $\pm 1$. On the other hand, if $\lfloor\frac{l}{k}\rfloor\equiv 1 \, (2)$, then the truncated block is of $\pm 1$ and we have $r=l-k\lfloor\frac{l}{k}\rfloor$ $\pm 1$. This concludes de proof.
    \end{proof}
\end{lema}

\begin{coro}\label{Coro:2ndCoro}
    For odd $k,l$ \begin{equation}\label{Eq:alpha_1Final}\alpha_1=\begin{cases}
        \frac{1}{l+k}\left(-\frac{1}{2}\left\lfloor\frac{l}{k}\right\rfloor -k\left(\left\lfloor\frac{l}{k}\right\rfloor +1 \right)+l-1\right) & \left\lfloor\frac{l}{k}\right\rfloor\equiv 0 \, (2)\\
        &\\
        \frac{1}{l+k}\left(-\frac{1}{2}\left(\left\lfloor\frac{l}{k}\right\rfloor -1\right)-l+k\left\lfloor\frac{l}{k}\right\rfloor-1\right)&\left\lfloor\frac{l}{k}\right\rfloor\equiv 1 \, (2)\\
    \end{cases}\end{equation}
    In particular, $\alpha_1\neq 0$ for any choice of coprime odd $k$ and $l$.

\begin{proof}
Note that $\alpha_1< 0$ since $-\frac{1}{2}\left\lfloor\frac{l}{k}\right\rfloor<0$, $k\left(\left\lfloor\frac{l}{k}\right\rfloor +1 \right)+l\leq 0$ and $-1<0$. Also, $-\frac{1}{2}\left(\left\lfloor\frac{l}{k}\right\rfloor -1\right)\leq 0$, $-l+k\left\lfloor\frac{l}{k}\right\rfloor<0$ and $-1<0$.
\end{proof}
\end{coro}

Combining Remark \ref{Lema:Converges_1/2}, Equation \eqref{Eq:alpha_1+1/2}, and Corollaries \ref{Coro:1stCoro} and \ref{Coro:2ndCoro} we obtain

\begin{teo}\label{Thm:Main2}
    Let $S=\{k,l\}\subset\mathbb{N}$ and let $\alpha_1$ be as in Equation \eqref{Eq:alpha_1Final}. The sequence \eqref{Eq:Main_kl} conveges (to $\frac{1}{2}$) only when $kl\equiv 0 \, (2)$. When $k,l\equiv 1 \, (2)$ the subsequences $b_n=a_{2n}$ and $c_n=a_{2n-1}$ converge to $\frac{1}{2}+\alpha_1$ and $\frac{1}{2}-\alpha_1$ respectively.
\end{teo}

\section{The General Case}

We now study the problem for a generic subtraction set $S=\{k_1,\ldots,k_t\}\subset\mathbb{N}$ with $k_i<k_j$ if $i<j$ and $\gcd\{k_1,\ldots,k_t\}=1$. In this situation, we shall prove that $(a_n^S)_{n\in\mathbb{N}_0}$ converges for every $S$ that has at least one even element (see Theorem \ref{Thm:MainGeneral}). However, we do not show that an all-odd subtraction set necesarily do not converge (although we conjeture it) nor we provide the exact limits of the even and odd subsequences for the (conjectured) non-convergent case.

Consider the general equation \eqref{Eq:Main}. We shall apply the same reasoning as in the previous section to identify when the sequences converge in terms of the properties of the roots of the characteristic polynomial $\chi_S(x)$ given in \eqref{Eq:Characteristic}. To simplify formulae, in the following we let $$d_i:=\begin{cases}k_t-k_i & 1\leq i\leq t-1 \\k_t & i=t\end{cases}.$$ Note that $\gcd\{k_1,\ldots,k_t\}=1$ if and only if $\gcd\{d_1,\ldots,d_t\}=1$.

We have the following (partial) generalization of Proposition \ref{Prop:Roots2}.

\begin{prop}\label{Prop:RootsGeneral}
 All the complex roots of $\chi_S(x)=tx^{d_t}+\sum_{i=1}^{t-1} x^{d_i}+1$ have norm less than 1, with the exception of the simple root $z=-1$ which happens only for the case $d_t\equiv 1 \, (2)$ and  $d_i\equiv 0 \, (2)$ for all $1\leq i\leq t-1$.
    \begin{proof}
        Let $z$ be a root of $\chi_S(x)$. Again by Lagrange's bound we have  \begin{equation}|z|\leq \max\left\{1,\sum_{i=0}^{d_t-1}\left|\frac{c_i}{c_{d_t}}\right|\right\}=1,\end{equation}
    where $c_0,\ldots,c_{d_t}$ denote the coefficients of $\chi_S(x)$. Assume $|z|=1$ and write $$tz^{d_t}+\sum_{i=1}^{t-1} z^{d_i}+1=0.$$ Now $tz^{d_t}$ lies in the circle of center $(0,0)$ and radius $t$ and $\sum_{i=1}^{t-1} z^{d_i}+1$ lies inside the circle of center $(1,0)$ and radius $t-1$. The only way that $\sum_{i=1}^{t-1} z^{d_i}+1$ can bring an element from distance $t$ to $0$ is if $z^{d_i}=1$ for all $i=1,\ldots,t-1$ and, consequently, $tz^{d_t}=-t$. We conclude that $z\in G_{d_i}$ for all $i=1,\ldots,t-1$ and $z\in G_{2d_t}$. That is, $$z\in \left(\bigcap_{i=1}^{t-1}G_{d_i}\right) \cap G_{2d_t}=G_{\gcd\{d_1,\ldots,d_{t-1},2d_t\}}=\begin{cases}
        G_1 & \exists\, d_i\equiv 1\, (2)\\
        G_2 & d_i\equiv 0\, (2)\, \forall\, i=1,\dots t-1
    \end{cases}$$
Now, $1$ is not a root of $\chi_S(x)$ so we conclude that $d_i\equiv 0\, (2)\, \forall\, i=1,\dots t-1$. Note that this implies that $d_t=k_t$ is odd, since if $k_t$ were even and also all the $d_i=k_t-k_i$ then all $k_i$ would be even, contradicting that $\gcd\{k_1,\ldots,k_t\}=1$. We conclude that $k_t$ must be odd and, hence, all $k_i$ with $i=1,\ldots,t-1$ must also be odd (for $d_i$ to be all even). We conclude that the only root of $\chi_S(x)$ of norm 1 is $z=-1$ for the case $k_i\equiv 1 \, (2)$ for all $i=1,\dots, t$. Furthermore, since $\chi'(-1)\neq 0$ then $-1$ is a simple root.
    \end{proof}
\end{prop}

\begin{teo}\label{Thm:MainGeneral}
    The sequence $a_n^S$ given in \eqref{Eq:Main} converges (to $\frac{1}{2}$) if $k_i\equiv 0\, (2)$ for some $i=1,\ldots,t$.
\begin{proof}
    For the mentioned case, Proposition \ref{Prop:RootsGeneral} asserts that all roots of $\chi_S(z)$ have norm less than $1$. Hence, $\lim_{n\rightarrow \infty}a_n=\frac{1}{2}$ by equation \eqref{Eq:ClosedForm}.
\end{proof}
\end{teo}

We now provide some remarks for the case not covered by Theorem \ref{Thm:MainGeneral}. On one hand, in order to determine if $(a_n^S)_{n\in\mathbb{N}_0}$ does or not converge in this setting it must be shown whether the coefficient $\alpha_1$ of $(-1)^n$ in Equation \eqref{Eq:ClosedForm} is or is not zero (respectively). We did this in the case $t=2$ by explicitly finding $\alpha_1$ and showing it is always negative. However, finding and explicity expression for this coefficient is rather difficult in this more general setting. We could proceed with the same arguments as in the case $t=2$ but under the assumption that all roots in this particular case are simple, which we believe is the case.
\begin{conj}\label{Conj:roots}
    The roots of $\chi_S(x)=tx^{d_t}+\sum_{i=1}^{t-1} x^{d_i}+1$ are all simple whenever $d_t\equiv 1 \, (2)$ and  $d_i\equiv 0 \, (2)$ for all $1\leq i\leq t-1$. 
\end{conj}
Indeed, we again solve the system \eqref{Eq:MatrixForm}, where $z_1=-1$. This system is invertible by Conjecture \ref{Conj:roots}, so the same arguments go through to compute the coefficients  $L_0(x)$ as in Equation \eqref{L_0}. In this case, an easy computation shows that  $$\frac{\chi_S(z)}{x+1}= \sum_{j=1}^t \sum_{i=d_{t-j+1}}^{d_{t-j}-1} j(-1)^i x^i.$$
Also, since $d_t$ is odd and $d_i$ are even for all $i=1,\ldots,t-1$, then $$\chi'(-1)=td_t(-1)^{d_t-1}+\sum_{i=1}^{t-1} (-1)^{d_i-1}d_i=td_t-\sum_{i=1}^{t-1} d_i=\sum_{i=1}^t d_i=\sum_{i=1}^t k_i.$$
Hence, since $\chi(z)=tf(x)$ we conclude that

$$L_0(x)=\dfrac{t}{k_1+\cdots+k_t}\left(\sum_{j=1}^t \sum_{i=d_{t-j+1}}^{d_{t-j}-1} \dfrac{j}{t}(-1)^i x^i\right).$$
So, the coefficients for $L_0(x)$ are:
$$L_{i0}=\begin{cases}
    \dfrac{j(-1)^i}{k_1+\cdots+k_t} & d_{t-j+1}\leq i< d_{t-j} \, (j=1,\ldots,t)
\end{cases}$$
We conclude that $$\alpha_1=\sum_{j=1}^t \sum_{i=d_{t-j+1}}^{d_{t-j}-1} \dfrac{j(-1)^i}{k_1+\cdots+k_t}\left(a_i^S-\frac{1}{2}\right).$$
This sum is challenging to address because the initial coefficients $a_i^S$ align with the first $k_t-1$ values of the sequence $a_n^{S-\{k_t\}}$, which no longer consists of blocks of zeros and ones. However, other non-necesarily explicit arguments could be found to detemine whether coefficient $\alpha_1$ in Equation \eqref{Eq:ClosedForm} is or is not null.

\begin{preg}\label{Question}
    Is $\alpha_1\neq 0$ for every $S=\{k_1,\ldots,k_t\}$ with $k_1<\cdots<k_t$ and $\gcd\{k_1,\ldots,k_t\}=1$?
\end{preg}

\section{Closing remarks}

In this section we outline some open problems,  potential avenues for further exploration, and conjectures that may suggest future research paths.

\subsection{The problem for all finite $S$.} In this paper we have addressed the study of the behaviour of the sequences $(a_n^S)_{n\in\mathbb{N}_0}$ for all finite $S$, completely determining it for $|S|< 3$ and partially for a generic $S$. A possible path to complete this last case would be to go through Conjecture \ref{Conj:roots} and Question \ref{Question} above (although probably no computing explicity the coefficient $\alpha_1$ corresponding to the root $-1$ of $\chi_S$). Additionally, for some particular choices of $S$ and for the complete case $|S|=3$, it may be possible even to find the limit of the even and odd subsequences of $(a_n^S)_{n\in\mathbb{N}_0}$.

\subsection{Dynamic Subtraction Sets}

A more general class of subtraction games is obtained by allowing the subtraction set depend on the oponent's previous move (see e.g. \cite{Schwenk}). For example, at a given position, the player whose turn it is to move can remove 1 chip or the entire pile; or a legal move is to remove any number of chips but not exceeding the amount removed by his opponent on the previous turn. Sequences derived from this problems would also satisfy a recurrence relation and can be attacked with some of the techniques presented in this article.

For example, the case with subtraction set $S=\{1,n\}$, which allows as a valid move to either remove 1 chip or the entire pile at each given position $n$, gives rise to the recurrence
$$
    a_n^S=\begin{cases}0&n=0 \\
    1-\frac{1}{2}a_{n-1}^S& n\geq 1 \end{cases}
$$
This is a linear recursion of order 1, so the solution can be computed straightforwardly as $\alpha_1z_1^n+\beta$ where $z_1$ is the root of the characteristic polynomial $\chi(x)=x-\frac{1}{2}$ and $\beta=\frac{2}{3}$ is a particular solution of the recursive formula. We find that $a_n^{S}=\alpha_1 (\frac{1}{2})^n+\frac{2}{3}$, which converges to $\frac{2}{3}$. This evidences that with dynamic substraction sets we no longer have converges to $\frac{1}{2}$.

On the other hand, the (trivial) case for $S=\{1,\ldots,n\}$ is straightforward. Indeed, in this situation:
$$
    a_n=\begin{cases}a_0=0,\, a_1=1& \\1-\dfrac{1}{n}(a_0+a_1+\cdots+a_{n-1}) & n\geq 2 \end{cases}
$$
It's easy to see by induction that $a_n=\frac{1}{2}$ for every $n\geq 2$. So the game is balanced for each $n\geq 2$.

\subsection{Several piles}
A next step can also be studying substraction games of more than one pile. In this case, and additional level of randomness must be considered as to which pile is selected to perform the move. Also, we may require different subtraction set for different piles. In the more general case, assume there are $r$ piles and we let $\mathcal{S}=\{S_j\}_{j=1,\dots,r}\subset\mathbb{N}$ be a collection of non-empty sets, set $S_j$ being the subtraction set for pile $j$. We represent a position of this game $(n_1,\ldots,n_r)$, pile $j$ containing $n_j$ chips. A similar argument as in Equation \eqref{Eq:Recursive} can be performed to obtain  multi-indexed sequences $(a_{n_1,\ldots,n_r}^{\mathcal{S}})_{n_i\in\mathbb{N}_0}$ representing the probability of winning for the next player to move from the position $(n_1,\ldots,n_r)$. Let $Prob(X,n_1,\ldots,n_r,j,i)$ stand for the conditional probability that player $X$ wins from position  $(n_1,\ldots,n_r)$ given that move $i$ in pile $j$ is selected. Then, for $n_i\geq \max\{S_i\}$, we have:
\begin{eqnarray}a_{n_1,\ldots,n_r}^{\mathcal{S}}&=&\displaystyle\sum_{j= 1}^r\sum_{i\in S_j} Prob(N,n_1,\ldots,n_{j-i},\ldots,n_r,j,i) \cdot Prob(\text{$i$ is selected}) \\
&=&\displaystyle\sum_{j= 1}^r\sum_{i\in S_j} (1-Prob(\text{$P$ wins from $(n_1,\ldots,n_{j-i},\ldots,n_r)$})) \cdot \frac{1}{r}\cdot \frac{1}{|S_j|}\nonumber  \\
&=&\displaystyle\sum_{j= 1}^r\sum_{i\in S_j} (1-a_{n_1,\ldots,n_{j-i},\ldots,n_r}^{\mathcal{S}})\cdot \frac{1}{r}\cdot \frac{1}{|S_j|}\nonumber\\
&=&\frac{1}{r}\displaystyle\sum_{j= 1}^r\frac{1}{|S_j|}\sum_{i\in S_j} (1-a_{n_1,\ldots,n_{j-i},\ldots,n_r}^{\mathcal{S}})\nonumber\\
&=&\frac{1}{r}\displaystyle\sum_{j= 1}^r \left(1-\frac{1}{|S_j|}\sum_{i\in S_j} a_{n_1,\ldots,n_{j-i},\ldots,n_r}^{\mathcal{S}}\right)\nonumber\\
&=&1-\frac{1}{r}\displaystyle \sum_{j= 1}^r \frac{1}{|S_j|}\sum_{i\in S_j} a_{n_1,\ldots,n_{j-i},\ldots,n_r}^{\mathcal{S}}\nonumber
\end{eqnarray}
It would be interesting to look into the stardard $r$=pile NIM, since in this case $S_j=\{1,\ldots,n\}$ for all $j=1,\ldots,r$ and the recursive relation becomes:

$$a_{n_1,\ldots,n_r}^{\mathcal{S}}=1-\frac{1}{rn}\displaystyle \sum_{j= 1}^r \sum_{i=1}^n a_{n_1,\ldots,n_{j-i},\ldots,n_r}^{\mathcal{S}}$$

\subsection*{Acknowledgements}The authors would like to thank Melina Privitelli for her useful contributions to the computational aspects of this paper.

\end{document}